\documentclass[11pt]{amsart}
\usepackage{graphicx}
\usepackage{amssymb}
\usepackage{amsmath}
\usepackage{amstext}

\usepackage{amsthm}
\usepackage{epsf}
\usepackage[scanall]{psfrag}
\usepackage{graphicx, color}

\def\Thmname{Theorem}
\def\Propname{Proposition}
\def\Lemmaname{Lemma}
\def\Definitionname{Definition}
\def\SSQ{
\left[%
\begin{array}{cc}
  e & f \\
  h & g \\
\end{array}%
\right] }

\newtheorem{Thm}{\Thmname}[section]

\newtheorem{Lem}[Thm]{\Lemmaname}

\def\onto{\twoheadrightarrow}

\def\J{\mathrel{{\mathcal J}}} 
\def\D{\mathrel{{\mathcal D}}} 

\def\Y{{\mathcal Y}}
\def\X{{\mathcal X}}

\def\R{\mathrel{{\mathcal R}}} 
\def\L{\mathrel{{\mathcal L}}} 
\def\H{\mathrel{{\mathcal H}}} 

\def\e<{\leq _{E}}

\font\petite=cmmi10 at 8pt
\def\malce{\mathbin{\hbox{$\bigcirc$\rlap{\kern-9pt\raise0,75pt\hbox{\petite
m}}}}}

\def\edge{\hskip-2pt\leftrightarrow\hskip-2pt}

\def\ctln{\centerline}

\def\ssk{\smallskip}
\def\msk{\medskip}
\def\bsk{\bigskip}

\def\veu{\vec{e}_1}
\def\vet{\vec{e}_2}

\def\vei{\vec{e}_i}
\def\vej{\vec{e}_j}
\def\vek{\vec{e}_k}
\def\ven{\vec{e}_n}
\def\hsk{\hskip5pt}

\title{Subgroups of free
idempotent generated semigroups: full linear monoids
 }

\author{Mark Brittenham}
\address{Department of Mathematics \\
University of Nebraska \\
Lincoln, Nebraska 68588-0323 \\
USA} \email{mbrittenham2@math.unl.edu}

\author{Stuart~W.~Margolis}
\address{Department of Mathematics \\
Bar Ilan University \\
52900 Ramat Gan \\
Israel} \email{margolis@math.biu.ac.il}

\author{John Meakin}
\address{Department of Mathematics \\
University of Nebraska \\
Lincoln, Nebraska 68588-0323 \\
USA} \email{jmeakin@math.unl.edu}

\thanks{The first author acknowledges support from NSF Grant
DMS-0306506. The second author acknowledges support from the
Department of Mathematics, University of Nebraska-Lincoln}

\begin{document}

\begin{abstract}

We develop some new topological tools to study maximal subgroups
of free idempotent generated semigroups. As an application, we
show that the rank 1 component of the free idempotent generated
semigroup of the biordered set of a full matrix monoid of size $n
\times n, n>2$ over a division ring $Q$ has maximal subgroup
isomorphic to the multiplicative subgroup of $Q$.

\end{abstract}

\maketitle

\section{Introduction}

Let $S$ be a semigroup with non-empty set $E = E(S)$ of
idempotents, and let $e, f \in E(S)$. It is easy to see that if $e
\in Sf \cup fS$, then both $ef$ and $fe$ must also be idempotents
of $S$. Products of idempotents of this form are referred to as
{\it basic products} in $E$. The set $E = E(S)$ relative to these
basic products forms a partial algebra that has been characterized
axiomatically as a {\it regular biordered set} by Nambooripad
\cite{Namb} in the case that $S$ is a (von-Neumann) regular
semigroup, and more generally as a {\it biordered set} by Easdown
\cite{Eas} for an arbitrary semigroup $S$. The basic products in
$E$ may be defined in terms of certain quasi-orders on $E$ that
are independent of the specific semigroup $S$ with biordered set
$E$. We refer to \cite{Namb} and \cite{Eas} for details of the
axiomatic characterization of biordered sets: we will not need
these details in the present paper.

Given any (axiomatically characterized) biordered set $E$, the
{\it free idempotent generated semigroup} on $E$ is the semigroup
$IG(E)$ with presentation

\begin{center}

$IG(E) = \langle E : e.f = ef$ if $ef$ is a basic product
$\rangle$.

\end{center}

Since a product of the form $e.e = e$ is clearly a basic product,
it is obvious that $IG(E)$ is an idempotent generated semigroup. A
theorem of Easdown \cite{Eas} shows that the biordered set of
idempotents of $IG(E)$ is $E$, that is, there is a bijection
between $E$ and the biordered set of idempotents in $IG(E)$ and
these two biordered sets have the same basic products. The
semigroup $IG(E)$ is a universal object in the category of
idempotent generated semigroups with biordered set $E$ and
morphisms that are one to one on idempotents. An analogous result
was proved earlier by Nambooripad \cite{Namb}, who constructed a
free regular idempotent generated semigroup $RIG(E)$ on a regular
biordered set $E$.

Given an idempotent $e$ of any semigroup $S$, the {\it maximal
subgroup} $H_{e}$ of $S$ with identity $e$ is the group of units
of the submonoid $eSe$ of $S$. For example, if $e$ is an
idempotent $n \times n$ matrix of rank $r $ in the monoid
$M_{n}(Q)$ of $n \times n$ matrices over a division ring $Q$, then
$H_{e} \cong GL_{r}(Q)$, the General Linear Group of size r over $Q$.
In \cite{BMM} it is shown that if $E$ is
a regular biordered set and $e \in E$, then the maximal subgroup
of $RIG(E)$ with identity $e$ is isomorphic to the maximal
subgroup of $IG(E)$ with identity $e$.  A question that has been
of interest in the literature is: which groups can arise as
maximal subgroups of a free idempotent-generated semigroup $IG(E)$
for some (axiomatically characterized) biordered set $E$? Early
results on this problem (for example \cite{Past2}, \cite{NP1})
suggested that such groups must be free, and in fact this
conjecture explicitly surfaced in the literature in \cite{McElw},
although it had been conjectured since the early 1980's.

In \cite{BMM}, the authors provided the first counterexample to
this conjecture, by showing that the free abelian group ${\Bbb
Z}\times{\Bbb Z}$ can arise in this context. Subsequently, Gray
and Ruskuc \cite{GR} have shown that {\it every} group arises in
this context. However, the structure of the maximal subgroups of
free idempotent generated semigroups on naturally occurring
biordered sets (such as the biordered set of the full linear
monoid $M_n(Q)$ over a division ring $Q$) is far from clear, and
the main purpose of the current paper is to provide some new
topological tools to study this problem.

We remark that idempotent generated semigroups arise naturally in
many parts of mathematics. First of all they are ``general". It is
well known that every semigroup $S$ embeds into a semigroup
generated by a set of idempotents of the same cardinality as $S$.
If $S$ is (finite) countable then $S$ embeds into a (finite)
semigroup generated by 3 idempotents. Idempotent generated
semigroups play an important part in the theory of reductive
algebraic monoids \cite{putcha, Renner}. Putcha's theory of
monoids of Lie type shows that one can consider the biordered set
of idempotents of such a monoid to be a generalization of a
building \cite{P1,P2} in the sense of Tits. Thus such objects
have a natural geometric structure.

In their paper \cite{BMM}, the authors defined the {\it
Graham-Houghton} $2$-complex $GH(E)$ of a regular biordered set
$E$, based on the work of Nambooripad \cite{Namb}, Graham
\cite{Gr} and Houghton \cite{Hough}, and they showed that the
maximal subgroups of $IG(E)$ are the fundamental groups of the
connected components of $GH(E)$. The $2$-cells of $GH(E)$
correspond to the {\it singular squares} of $E$ defined by
Nambooripad \cite{Namb}. Gray and Ruskuc \cite{GR} give a new
proof that singular squares give presentations of maximal
subgroups of $IG(E)$ for an arbitrary biordered set $E$.

We outline the main topological idea of this paper. Let $S$ be an
idempotent generated regular semigroup with biordered set
$E=E(S)$. Then there is a surjective idempotent separating
morphism $f:RIG(E) \onto S$ that is an isomorphism on $E$. It
follows that for every maximal subgroup $G$ of $S$ there is a
unique connected component $C$ of $GH(E)$ and a unique point $x$ of $C$ and a
surjective morphism $\pi_{1}(C,x) \onto G$. By basic algebraic
topology, $G$ acts transitively and fixed point free on each fibre
of the connected cover $C(G)$ of $C$ that has fundamental group
the kernel of this latter morphism. Thus if $C(G)$ is simply
connected, $G$ is isomorphic to a maximal subgroup of $RIG(E)$.

We apply these ideas to show that if $Q$ is a division ring, then
the maximal subgroup of $IG(E(M_{n}(Q)))$ corresponding to an
idempotent matrix of rank $1$ is $Q^*$, the multiplicative group
of units of $Q$. The maximal subgroup of $IG(E(M_{n}(Q)))$
corresponding to an idempotent matrix of rank $n-1$ is a free
group, but the structure of the maximal subgroups of
$IG(E(M_{n}(Q)))$ corresponding to idempotent matrices of rank $k$
for $1 < k < n-1$ remains far from clear.

\section{Preliminaries}

Recall that a semigroup $S$ is called {\it regular} if $a \in aSa$
for each $a \in S$. In a very influential paper \cite{Namb},
Nambooripad studied the structure of regular semigroups via his
theory of {\it inductive groupoids}. He found an axiomatic
characterization of the set of idempotents of a regular semigroup
$S$ relative to the basic products in $E(S)$ as a ``regular
biordered set" and he described the inductive groupoid associated
with the free regular idempotent generated set $RIG(E)$ on a
regular biordered set $E$.  A presentation for $RIG(E)$ was
provided by Pastijn \cite{Past2}. In \cite{BMM}, the authors
showed that if $E$ is a regular biordered set, then the maximal
subgroups of $RIG(E)$ are isomorphic to the maximal subgroups of
the semigroup $IG(E)$ defined above. We refer the reader to
\cite{Namb} and \cite{BMM} for details. {\it We shall assume
throughout the remainder of this paper that all biordered sets
under consideration are biordered sets of regular semigroups}.

Recall that the Green's relations $\R$ and $\L$ on a semigroup $S$
are defined by $a \R b$ iff $aS^1 = bS^1$ and $a \L b$ iff $S^{1}a
= S^{1}b$. When restricted to $E = E(S)$ these are defined by
basic products: $e \R f$ iff $ef = f$ and $fe = e$, and $e \L f$
iff $ef = e$ and $fe = f$. Thus, one can consider the $\R$ and
$\L$ relations on an arbitrary (axiomatically defined) biordered
sets. By the theorems of Nambooripad and Easdown mentioned above,
these are exactly the restrictions of the corresponding Green's
relations on $RIG(E)$ and $IG(E)$ to $E$.

We can define Green's relation $\D$ on a biordered set $E$ as the
transitive closure of $\R \cup \L$. It follows from the work of
Fitz-Gerald \cite{FitzG} and Nambooripad \cite{Namb} that this is
the restriction of Green's relation $\D$ to the idempotents of
$RIG(E)$. Thus we identify these relations on $E$ with the
classical notions on $RIG(E)$ and $IG(E)$ without further mention.
In Nambooripad's language, \cite{Namb}, two idempotents of
$RIG(E)$ or $IG(E)$ are $\D$ related if and only if there is an
$E$-path between them. This is just another way of saying that
$\D$ is the transitive closure of $\R \cup \L$ restricted to $E$.

Central to Nambooripad's construction of the inductive groupoid of
$RIG(E)$ is the notion of a {\it singular square} in the (regular)
biordered set $E$. An {\it $E$-square} is a sequence $(e,f,g,h,e)$
of elements of $E$ with $e\R f \L g \R h \L e$. Unless otherwise
specified, we will assume that all $E$-squares are non-degenerate,
i.e. the elements $e,f,g,h$ are all distinct. An idempotent $t=t^2
\in E$ {\em left to right singularizes} the $E$-square
$(e,f,g,h,e)$ if $te=e, th=h, et=f$ and $ht=g$ where all of these
products are defined in the biordered set $E$. Right to left, top
to bottom and bottom to top singularization is defined similarly
and we call the $E$-square {\em singular} if it has  a
singularizing idempotent of one of these types. All of these
products are basic products, so they make sense in any semigroup
with biorder isomorphic to $E$.

The following simple but important fact was first noted by
Nambooripad \cite{Namb}. Recall that the right zero semigroup on a
set $X$ is the set $X$ with multiplication $xy=y$ for all $x,y \in
X$. A left zero semigroup is the dual notion and a rectangular
band is the direct product of some left zero semigroup with some
right zero semigroup. Thus a rectangular band is defined on some
set of the form $X \times Y$ with multiplication $(x,y)(x\prime ,
y\prime) = (x,y\prime)$, for all $x,x\prime \in X, y,y\prime \in
Y$.

\begin{Lem} \label{Trivsing}
Let $(e,f,g,h,e)$ be a singular $E$-square in a semigroup $S$.
Then $efghe=e$: in other words, $\SSQ$ is a rectangular band in
any semigroup with biordered set $E$.

\end{Lem}

We remark that the converse of Lemma \ref{Trivsing} is obviously
false. For example if the semigroup $S$ is just a $2 \times 2$
rectangular band $\SSQ$, then it is not a singular square since
there is no idempotent available to singularize it.

Recall from \cite{BMM} that the {\it Graham-Houghton} graph of a
(regular) biordered set $E$ is the bipartite graph  with vertices
the disjoint union of the set of $\R$-classes of $E$ and the set
of $\L$-classes of $E$, and with a directed (positively oriented)
edge from an $\L$-class $L$ to an $\R$-class $R$ if there is an
idempotent $e \in L \cap R$ (and a corresponding inverse edge from
$R$ to $L$ in this case). We now add 2-cells to this graph, one
for each singular square $(e,f,g,h,e)$. Given this square we sew a
2-cell onto this graph with boundary $ef^{-1}gh^{-1}$. We call the
resulting 2-complex the {\it Graham-Houghton complex} of $E$ and
denote it by $GH(E)$.

The following theorem in \cite{BMM} is based on the work of
Nambooripad, and is the principal tool used in \cite{BMM} to
construct maximal subgroups of free idempotent-generated
semigroups on biordered sets.

\begin{Thm} \cite{BMM} \label{BMM} Let $E$ be a regular biordered set.
Then the maximal subgroup of $IG(E)$ based at $e \in E$ is
isomorphic to the fundamental group $\pi_{1}(GH (E), L_{e})$ of
the Graham-Houghton complex of $E$ based at $L_{e}$.
\end{Thm}

The following theorem is used crucially in this paper.

\begin{Thm} \label{ConnComp}

Let $E$ be a regular biordered set. Then for any regular
idempotent generated semigroup $S$ with $E \approx E(S)$, there is
a one to one correspondence between connected components of
$GH(E)$ and $\D$ classes of $S$. Every maximal subgroup $G$ of $S$
is a quotient of a unique maximal subgroup of $RIG(E)$ that
belongs to the unique connected component corresponding to the
$\D$-class of $G$ of $S$.

\end{Thm}

{\bf Proof} Since $E \approx E(S)$, there is a surjective
idempotent separating morphism $\varphi:RIG(E) \onto S$. It is well
known that $f$ then induces a bijection between $RIG(E)/
\mathcal{K}$ to $S/ \mathcal{K}$ for any of Green's relations
$\mathcal{K}$. In particular, this is true for $\mathcal{K}= \D$.
As mentioned above, two idempotents are $\D$ related in $RIG(E)$
if and only if there is an $E$-path between them in $E$. Since positive
edges of $GH(E)$ are in one to one correspondence with $E$, a
straightforward induction on the length of a path shows that there
is a path between two vertices of $GH(E)$ if and only if there is
an $E$-path in $E$ whose first edge belongs to the first vertex
(recall that the vertices of $GH(E)$ are the disjoint union of
$\R$ and $\L$ classes of $E$) and whose last edge belongs to the
last vertex. Since idempotent separating morphisms between regular
semigroups induce a 1-1 correspondence between $\H$ classes, the
second statement of the theorem follows. This completes the proof.

\section{A Freeness Criterion}

Let $S$ be a regular idempotent generated semigroup and let
$E=E(S)$ be its biordered set of idempotents. Let $G$ be a maximal
subgroup of $S$. By Theorem \ref{BMM} and Theorem \ref{ConnComp}
there is a unique connected component $\mathcal{G}$ of $GH(E)$,
and a surjective morphism $f:\pi_{1}(\mathcal{G}) \onto G$ from
the fundamental group of $\mathcal{G}$ to $G$.

Therefore, there is a unique (up to isomorphism) connected cover
$\mathcal{C}(G)$ of $\mathcal{G}$ that has $Ker(f)$ as
fundamental group. It is well known, \cite{hatch} that $G$ acts
freely on the fibres of $\mathcal{C}(G)$ and that
$\mathcal{G} \approx \mathcal{C}(G)/G$. The construction
of $\mathcal{C}(G)$ is a simple exercise in covering
theory, but since we need the details, we include them here.

Let $\mathcal{C}$ be a connected 2 complex and let
$\varphi:\pi_{1}(\mathcal{C}) \onto G$ be a surjective morphism to a
group $G$. Let $\mathcal{T}$ be a spanning tree of
$V(\mathcal{C})$. Then $\pi{_{1}}(\mathcal{C})$ and $G$ are
generated by the positive edges (relative to some orientation) not
in $\mathcal{T}$.

For convenience if $e$ is an edge of $\mathcal{C}$ we identify $e$
with its value as a generator in $\pi_{1}(\mathcal{C})$ as above
and $\varphi(e)$ the corresponding value in $G$. Note that $e$ and
$\varphi(e)$ are the identity element if $e \in \mathcal{T}$

Define the cover $\mathcal{C}(G)$ as follows:

{\bf Vertices:} $G \times V(\mathcal{C})$

{\bf Edges:} We have an edge from $(g,x)$ to $(h,y), g,h \in G,
x,y \in V(\mathcal{C})$ iff there is an edge $e$ from $x$ to $y$
in $\mathcal{C}$ and $g\cdot\varphi(e)=h$ in $G$.

{\bf 2-Cells:} For every cell of $\mathcal{C}$ sewed on as a loop
in $\mathcal{C}$ from some vertex $x$ to itself, we sew a copy of
it as a loop from $(g,x)$ to itself following the edges modified
by their lifts in the definition of edges above.

Then this is a cover of $\mathcal{C}$ denoted by $\mathcal{C}(G)$
under the projection $G \times V(\mathcal{C}) \rightarrow
V(\mathcal{C})$. $G$ acts transitively and fixed point free on the
left of each fibre by $g(h,x)=(gh,x)$ and the quotient
$\mathcal{C(G)}/G \approx \mathcal{C},
\pi{_{1}}{(\mathcal{C}(G))}\unlhd \pi{_{1}}(\mathcal{C})$ and $G
\approx \pi{_{1}}(\mathcal{C})/\pi{_{1}}(\mathcal{C}(G))$

Finally, $\mathcal{C}(G)$ is clearly the universal object in the
category of all covers of $\mathcal{C}$ on which $G$ acts
transitively on each fibre.

{\bf Example 1} If $\mathcal{C}$ is a bouquet of $X$ circles (no 2
cells), then $G$ is an $X$ generated group and $\mathcal{C}(G)$ is
just the Cayley graph of $G$ relative to the presentation
$f:\pi{_{1}}(\mathcal{C}) \rightarrow G$.

{\bf Example 2} Let $E$ be a regular biordered set and let $S$
be a regular idempotent generated semigroup with $E(S)\approx E$.
Let $G$ be a maximal subgroup of $S$. Then $G$ corresponds by
Theorem \ref{ConnComp} to a unique connected component  $\mathcal{C}$ of $GH(E)$
and there is a surjective morphism from the fundamental group
corresponding to this component to $G$.

Putting all this together we get the following criterion for $G$ to be isomorphic to a 
maximal subgroup of $RIG(E)$ (or equivalently $IG(E)$ \cite{BMM}). We use
the notation in Example 2.

{\bf Freeness Criterion} If the cover $\mathcal{C}(G)$ of the
group $G$ is simply connected, then $G$ is isomorphic to the maximal subgroup of
the corresponding component of $RIG(E)$.

\section{Matrices over division rings}

Throughout this section, $Q$ will be a division ring, $M_{n}(Q)$
will denote the full linear monoid of $n \times n$ matrices over
$Q$ and $GL_n(Q)$ will denote the general linear group, i.e.  the
group of units of $M_n(Q)$. We will use both lower case and upper
case letters to denote matrices.
We make use of the covering space methods in the previous section
to study the maximal subgroups of the free idempotent generated
semigroup on the biordered set of idempotents of this monoid. In
particular, we prove the following theorem, which is the main
result in this section.

\begin{Thm} \label{rankone}

Let $E$ be the biordered set of $M_{n}(Q)$, for $Q$ a division
ring, and let $e$ be an idempotent matrix of rank $1$ in
$M_{n}(Q)$. For $n \geq 3$, the maximal subgroup of $IG(E)$ with
identity $e$ is isomorphic to $Q^*$, the multiplicative group of
units of $Q$.

\end{Thm}

For basic facts about matrices over division rings, see the book
by Jacobson \cite{jacobs}.
There is a great deal of information about full linear monoids,
particularly in the case where $Q$ is a field (see for example the
books of Putcha \cite{putcha} and Okni\'nski \cite{okn}).
Much of the basic structural information about full linear monoids
over fields extends to the case where $Q$ is a division ring, but
care must be taken to extend some of these results if $Q$ is not
commutative. For example, linear combinations of rows will always
be considered using left scalar multiplication, and linear
combinations of columns using right scalar multiplication;
consequently row spaces are left row spaces and column spaces are
right column spaces.

It is well known that the set of matrices of a fixed rank $k \leq
n$ forms a $\J$-class in the monoid $M_{n}(Q)$. Here the rank of a
matrix $a \in M_{n}(Q)$ is the (left) row rank of $a$, which is
the same as the (right) column rank of $a$. In fact, for matrices
$a,b \in M_{n}(Q)$, we have $a \J b$ iff $GL_{n}(Q) \, a \,
GL_{n}(Q) = GL_{n}(Q) \, b \, GL_{n}(Q)$ iff $rank(a) = rank(b)$,
and $\J = \D$. The maximal subgroup of the $\J$-class of all
matrices of rank $k$ is isomorphic to $GL_{k}(Q)$.  The Green's
$\R$ and $\L$ relations on $M_{n}(Q)$ are characterized by

 $a \R b$ iff $a \, GL_{n}(Q) = b \, GL_{n}(Q)$ iff
 $Col(a) = Col(b)$, and

 $a \L b$ iff $GL_{n}(Q) \, a = GL_{n}(Q) \, b$  iff
$Row(a) = Row(b)$.

Let $D_k$ be the $\D$-class of $M_n(Q)$ consisting of the rank $k$
matrices, and $D_{k}^0$ the corresponding completely 0-simple
semigroup. Let $\Y_k$ be the set of all matrices of rank $k$ which
are in reduced row echelon form and let $\X_k$ be the set of
transposes of these matrices.
The structure of $D_{k}^0$ is described in the following theorem
(see \cite{okn}).

\begin{Thm} \label{rees}
$D_{k}^0 \cong {\mathcal M}^{0}(\X_k,GL_{k}(Q),\Y_k,C_k)$ where
the matrix $C_{k} = (C_k(y,x))$ is defined for $x \in \X_{k}, y\in
\Y_k$ by $C_k(y,x) = yx$ if $yx$ is of rank $k$ and $0$ otherwise.

\end{Thm}

By Theorem \ref{rees} and the basic structure of Rees matrix
semigroups (see, for example \cite{CP}), every matrix $a$ of rank
$k$ can be uniquely expressed in the form $a = xhy$ where $x \in
\X, y \in \Y$ and $h$ is a block diagonal matrix of the form $
\left[%
\begin{array}{cc}
  h' & 0 \\
  0 & 0 \\
\end{array}%
\right] $ where $h' \in GL_k(Q)$.
Since $x$ has $n-k$ columns of zeroes at the right of the matrix
and $y$ is the transpose of a matrix of this form we see that we
may write the matrix $a = xhy$ of rank $k$ in the form $a = vw^T$
for some $n \times k$ matrices $v,w$ of rank $k$ (choose $v =
x'h'$ where $x'$ is obtained from $x$ by deleting the last $n-k$
columns, and $w^T$ is obtained from $y$ by deleting the last $n-k$
rows). Also, we may replace $x$ [resp. $y$] in the above by
matrices of the form $x_{1} = xh_1$ [resp. $y_{2} = h_{2}y$] for
any matrices $h_{1},h_{2}$ in the maximal subgroup of $
\left[%
\begin{array}{cc}
  I_{k} & 0 \\
  0 & 0 \\
\end{array}%
\right] $ and get isomorphic Rees matrix semigroups. Thus we may
replace $yx$ by $w^Tv$ as above in the definition of the matrix
$C_k$ and obtain an isomorphic Rees matrix semigroup.

If $a$ is a rank-$k$ matrix, expressed as $a=vw^T$ as above, then
it is routine to check that  $a$ is an idempotent  iff $w^Tv=I_k$.
Two rank-$k$ $n\times n$ matrices $a=v_1w_1^T,b=v_2w_2^T$ are
 $\L$-related iff $a$ and $b$ have the same row space, which in
turn is true iff $w_1^T=mw_2^T$ for some non-singular $k\times k$
matrix $m$. Similarly, $a$ and $b$ are $\R$-related iff they have
the same column space, i.e., $v_1=v_2m$ for some non-singular
$k\times k$ matrix $m$. An $\L$-class can therefore be identified
with the equivalence class $[w^T]=\{mw^T : m\in GL_k(Q)\}$, and an
$\R$-class can be identified with $[v] = \{vm : m\in GL_k(Q)\}$
where $v,w$ are $n \times k$ matrices of rank $k$.

Rectangular bands in the biordered set of idempotents $E =
E(M_{n}(Q))$ may be characterized from the representation of rank
$k$ matrices described above.

If $e\in [v_1]\cap[w_1^T],f\in [v_2]\cap[w_1^T],g\in
[v_2]\cap[w_2^T],h\in[v_1]\cap[w_2^T]$ forms an $E$-square, then
$e=v_1(w_1^Tv_1)^{-1}w_1^T,f=v_2(w_1^Tv_2)^{-1}w_1^T,
g=v_2(w_2^Tv_2)^{-1}w_2^T,h=v_1(w_2^Tv_1)^{-1}w_2^T$.  Then
this $E$-square is a rectangular band iff $efghe = e$. A
calculation of this product shows that this happens iff

$$ (*) \ldots (w_1^Tv_2)(w_2^Tv_2)^{-1}(w_2^Tv_1)(w_1^Tv_1)^{-1}=I_k.$$

(Note that the identity (*) is independent of the choice of
representatives of $v_1,v_2,w_1,w_2$ in their equivalence
classes.) We have the following rather pleasant fact about the
semigroup $M_{n}(Q)$.

\begin{Thm} \label{rectsing}

Every non-trivial rectangular band in $M_{n}(Q)$ (for a division
ring $Q$) is a singular square.

\end{Thm}

{\bf Proof} Given an $E$-square
$\left[%
\begin{array}{cc}
  e & f \\
  h & g \\
\end{array}%
\right]$ consisting of  $n\times n$ idempotent matrices of rank
$k$ with coefficients in the division ring $Q$, we can, by
conjugating $e$ by a change of basis matrix whose columns are a
basis for the column space of $e$ followed by a basis for the
nullspace of $e$, assume that
$e=\left[%
\begin{array}{cc}
  I & 0 \\
  0 & 0 \\
\end{array}%
\right]$, where $I$ = the $k\times k$ identity matrix and the
$0$'s are matrices of $0$'s of the appropriate size. Then since
$ef=f,fe=e$, etc., a routine calculation demonstrates that after
the same conjugation
we have $f=\left[%
\begin{array}{cc}
  I & b \\
  0 & 0 \\
\end{array}%
\right]$, $h=\left[%
\begin{array}{cc}
  I & 0 \\
  a & 0 \\
\end{array}%
\right]$, and $g=\left[%
\begin{array}{cc}
  I & b \\
  a & ab \\
\end{array}%
\right]$. Finally, computing $g^2=g$ we find that $I+ba=I$ (in the
upper left corner), so $ba=0$.

We now show that for this quartet of idempotents there is a matrix
$\eta$ which singularizes the $E$-square. Conjugating this matrix
by the inverse of our change of basis matrix gives a matrix
singularizing the original $E$-square.
Writing $\eta=\left[%
\begin{array}{cc}
  x & y \\
  z & c \\
\end{array}%
\right]$ and computing the required (for left-to-right
singularization) products $\eta e=e$, $e\eta=f$, etc., together
with $\eta^2=\eta$,
we find that $\eta$ must have the form $\eta=\left[%
\begin{array}{cc}
  I & b \\
  0 & c \\
\end{array}%
\right]$ with $bc=0$, $ca=a$, and $c^2=c$; further, any such
matrix will be an idempotent singularizing the $E$-square. We now
proceed to construct the needed matrix $c$.

Recalling that column spaces always refer to right-linear
combinations of columns, note first that the condition $ba=0$ is
equivalent to $b($col$(a))=\{0\}$, which is equivalent to
col$(a)\subseteq$ null$(b)$. Similarly, the requirement $bc=0$ therefore
requires col$(c)\subseteq$ null$(b)$.

The condition $ca=a$ implies that the columns of $a$ are
right-linear combinations of the columns of $c$, so
col$(a)\subseteq$ col$(c)$. But conversely, since $c$ is
idempotent, col$(a)\subseteq$ col$(c)$ implies that $ca=a$. To see
this, since every column of $a$ is a right-linear combination of
the columns of $c$ we have $a=cs$ for some $n\times n$ matrix $s$.
Then $ca=c(cs)=(c^2)s=cs=a$, as desired.

So our requirements for the matrix $c$ are: $c$ is idempotent
($c^2=c$) and col$(a)\subseteq$ col$(c)\subseteq$ null$ (b)$.
Since every ${\R}$-class of $M_n(Q)$ contains an idempotent we can
arrange to have col$(a) = $col$(c)$ (and then col$(c) =
$col$(a)\subseteq$ null$(b)$ is immediate).

\bigskip

This result enables us to complete the description of the
Graham-Houghton complex $K=GH(E)$ of the set of idempotents $E$ of
the semigroup $M_{n}(Q)$ of $n\times n$ matrices over the division
ring $Q$. The vertices consist of the equivalence classes
$[v]$,$[w^T]$ of sets of $k$ linearly independent column (resp.,
row) vectors, $0\leq k\leq n$; $[v_1]=[v_2]$ iff $v_1$ and $v_2$
have the same column space, i.e., $v_2=v_1x$ for some non-singular
$k\times k$ matrix $x$; $[w_1^T]=[w_2^T]$ iff they have the same
row space, i.e., $w_2^T=xw_1^T$ for $x\in GL_k(Q)$. There is an
edge joining $[v]$ and $[w^T]$ iff $w^Tv$ is a non-singular
$k\times k$ matrix. A non-degenerate 4-cycle
$([v_1],[w_1^T],[v_2],[w_2^T])$ (so
$[v_1]\neq[v_2],[w_1^T]\neq[w_2^T]$) bounds a 2-cell iff
$(w_1^Tv_2)(w_2^Tv_2)^{-1}(w_2^Tv_1)(w_1^Tv_1)^{-1}=I_k$. These
constitute all of the edges and 2-cells in the complex.

\medskip

Using this same notation, we can now describe the cover $\widetilde{K}_{n,k}$, defined in section 3 (and
referred to as $\mathcal{C}(G)$ there)
corresponding to the $\D$-class of the rank-$k$ matrices.
$\widetilde{K}_{n,k}$ has
vertices the pairs $(g,[w^T])$,$([v],h)$, where $g,h\in GL_k(Q)$;
since column spaces are right vector spaces, $G$ acts on $v$ on the
right, and so we will write the latter pairs as $([v],h)$.
If for each equivalence class we choose a (fixed) representative
$w_0^T,v_0$ then we can identify $(g,[w^T])=(g,[w_0^T])$ with
$gw_0^T$ and $([v],h)$ with $v_0h^{-1}$.
As $g,h$ range over
$GL_k(Q)$, this identifies the vertices of $\widetilde{K}_{n,k}$
with the set of all rank-$k$ $k\times n$ and $n\times k$ matrices,
respectively.
In the notation of section 3, our morphism
$\varphi:\pi_1(GH(E))\rightarrow GL_k(Q)$ is $([v],[w^T])=w_0^Tv_0$,
where $([v],[w^T])$ is the edge from $[v]$ to $[w^T]$.
There is an edge from $gw_0^T$ to $v_0h^{-1}$ iff
$g(w_0^Tv_0)=h$, that is, $(gw_0^T)(v_0h^{-1})=I_k$, where $I_k$
is the $k\times k$ identity matrix. So the vertices of
$\widetilde{K}_{n,k}$ consist of the rank-$k$ $k\times n$ and
$n\times k$ matrices, and there is an edge from  $w^T$  to $v$ iff
$w^Tv=I_k$. Finally, there is a 2-cell with boundary any 4-cycle
in the 1-skeleton of $\widetilde{K}_{n,k}$.

\bigskip

{\bf Proof of Theorem \ref{rankone}.} We denote by $K_{n,1}$ the
subcomplex of $K$ spanned by the rank-$1$ vertices. By Theorem
\ref{BMM} and the Freeness Criterion of section 3, we will be able to
prove Theorem \ref{rankone}, if we can show that the cover
$\widetilde{K}_{n,1}$ is simply connected. By construction,
$\widetilde{K}_{n,1}$ has vertex set consisting of all nonzero
$n\times 1$ (column) vectors $v$ and all nonzero $1\times n$ (row)
vectors $w^T$. There is an edge $v\leftrightarrow w^T$ (consisting
of a positively oriented edge from  $w^T$ to $v$ and its inverse
edge from $v$ to $w^T$ ) iff $w^Tv=1$. Finally, each 4-cycle in
$\widetilde{K}_{n,1}^{(1)}$ is the boundary of a 2-cell in
$\widetilde{K}_{n,1}$.

\bigskip

To show that $\widetilde{K}_{n,1}$ is simply connected, that is,
that $\pi_1(\widetilde{K}_{n,1})=\{1\}$, we need to show that
every loop in $\widetilde{K}_{n,1}$ is null-homotopic. More
precisely, choosing a maximal tree $T$ in
$\widetilde{K}_{n,1}^{(1)}$, $\pi_1(\widetilde{K}_{n,1})$ is
generated by loops, one for each edge $\epsilon$ not in $T$. The
loops start at the basepoint, run out the tree to one endpoint of
$\epsilon$, across $\epsilon$, and then back in the tree to the
basepoint. It suffices to show that each of these loops is
null-homotopic, and for this it is enough to show that each edge
$\epsilon$ is homotopic, rel endpoints, to an edge path in $T$. It
is this last statement which we will now prove. We will
carry out this verification in steps, in the process building the
tree $T$ in steps as well.

The basic shortcut which we will use is the following observation.
If $T\subseteq \widetilde{K}_{n,1}^{(1)}$ is a tree and
$\epsilon_1,\ldots ,\epsilon_n,\epsilon\in K^{(1)}$ are edges with
all endpoints lying in $T$, and if each $\epsilon_i$ is homotopic
in $K$, rel endpoints, to an edgepath in $T$, \underbar{and}
$\epsilon$ is homotopic in $K$, rel endpoints, to an edgepath
$\gamma$ in $T\cup \epsilon_1\cup\cdots\cup \epsilon_n$, then
$\epsilon$ is also homotopic, rel endpoints, to an edgepath
$\delta$ in $T$. This is because we can concatenate a sequence of
homotopies, each supported on an edge $\epsilon_i$ lying in the
edgepath $\gamma$, deforming $\epsilon_i$ into $T$, to further
homotope $\epsilon$ into $T$. The net effect of this observation
is that, in the course of our proof, anytime we have shown that an
edge $\epsilon$ can be deformed into our tree $T$, we can act as
if $\epsilon$ were actually \underbar{in} $T$ and build our
further deformations to map into the union of $T$ and $\epsilon$
(and all other edges we have shown can deform into $T$). To
reinforce this, we will talk of the edges of our tree as being
colored ``green'', and say that any edge that we can deform into
$T$ has turned green. Then to continue to move our proof forward
we are required only to show that any further edge can be deformed
into the green edges. In this way more and more edges become
green; the proof ends when we have shown that every edge can be
turned green.

\msk

Our approach will be to choose nested collections ${\mathcal V}_i$
of vertices, and then, inductively, extend the tree $T_{i-1}$ from
a tree with the previous vertex set ${\mathcal V}_{i-1}$ to a tree
whose vertex set is ${\mathcal V}_i$, and show that the edges of
the full subcomplex of $\widetilde{K}_{n,1}$ with vertex set
${\mathcal V}_i$ can all be turned green. These edges then can
automatically be assumed to be green when moving on to the next
vertex set ${\mathcal V}_{i+1}$; a homotopy rel endpoints into the
tree $T_i$ is also a homotopy into the tree $T_{i+1}$. Our
approach relies on the fact that the homotopies take place across
the 2-cells of $\widetilde{K}_{n,1}$, whose boundaries are the
4-cycles in $\widetilde{K}_{n,1}^{(1)}$. Any time that we can find
a 4-cycle three of whose sides have turned green, the square then
provides a homotopy of the fourth side into the green edges,
enabling us to turn the fourth side green, as well. In what
follows, we will signify this by stating that the sequence of
edges $v_1\edge v_2\edge v_3\edge v_4$ ``yields'' the edge
$v_1\edge v_4$, meaning that it can now be turned green. Our proof
essentially consists of finding a way to list the edges in each of
our full subcomplexes so that, for each edge $\epsilon$ outside of
the tree $T_i$, there is a square containing $\epsilon$ so that
the other three edges of the square are each either in the tree
$T_i$ or appear earlier in the list, and so, by induction, can be
assumed to have turned green. This enables us to turn the edge
$\epsilon$ green, as well, and continue the induction.

In what follows we denote by $\vec{e}_i$ the vector with $1$ in
the $i$-th coordinate and $0$ in the remaining coordinates.

\msk

We start with the case $n=3$; our last induction will be on $n$.
Our first vertex set ${\mathcal V}_1$ consists of the nonzero
vectors $v,w^T$ all of whose entries are 0 or 1, with at least one
0. We build the tree $T_1$ by adding the edge from every $w^T$
with first coordinate 1 to $\veu$ and then an edge from every $v$
with first coordinate 1 to $\veu^T$, and then add, for every
vertex $v,w^T$ whose first non-0 coordinate occurs in the
$i^{\text{th}}$ entry, $i>1$, the edge from $v,w^T$ to $\veu+\vei$
and $\veu^T+\vei^T$, respectively. (Note that this has,
implicitly, already used the hypothesis $n\geq 3$, so that
$\veu+\vei\in{\mathcal V}_1$.) Since, for every edge in $T_1$, at
the time it is added exactly one of its (non-$\veu$) endpoints
does not yet lie in the part of $T_1$ constructed up to that
point, their union forms a tree, by induction. The edges already
in $T_1$ are therefore

\ssk

$\left(1\ 0\ 0\right)\edge$ each of $\left(1\ 0\
0\right)^T,\left(1\ 1\ 0\right)^T, \left(1\ 0\ 1\right)^T$,

$\left(1\ 0\ 0\right)^T\edge$ each of $\left(1\ 1\ 0\right),
\left(1\ 0\ 1 \right)$,

$\left(1\ 1\ 0\right)^T\edge$ each of $\left( 0\ 1\
0\right),\left(0\ 1\ 1\right)$ and

$\left(1\ 1\ 0\right)\edge$ each of $\left( 0\ 1\
0\right)^T,\left(0\ 1\ 1\right)^T$,

$\left(1\ 0\ 1 \right)\edge\left(0\ 0\ 1\right)^T$ and $\left(1\
0\ 1 \right)^T\edge\left(0\ 0\ 1\right)$.

\ssk

\noindent By inspection, the remaining edges joining vertices in
${\mathcal V}_1$ are

\ssk

$\left(1\ 1\ 0\right)\edge\left(1\ 0\ 1\right)^T$ and $\left(1\ 1\
0\right)^T\edge\left(1\ 0\ 1\right)$,

$\left(1\ 0\ 1\right)\edge\left(0\ 1\ 1\right)^T$ and $\left(1\ 0\
1\right)^T\edge\left(0\ 1\ 1\right)$,

$\left(0\ 1\ 0\right)\edge\left(0\ 1\ 0\right)^T$,$\left(0\ 1\
1\right)^T$ and $\left(0\ 1\ 0\right)^T\edge$,$\left(0\ 1\
1\right)$,

$\left(0\ 0\ 1 \right)\edge\left(0\ 1\ 1\right)^T$,$\left(0\ 0\
1\right)^T$ and $\left(0\ 0\ 1 \right)^T\edge\left(0\ 1\
1\right)$.

\ssk

\noindent Then the following 4-cycles show how to turn each of
these edges, in turn, green:

\ssk

$\left(1\ 1\ 0\right)\edge \left(1\ 0\ 0\right)^T\edge \left(1\ 0\
0\right)\edge \left(1\ 0\ 1\right)^T$ \hsk yields \hsk $\left(1\
1\ 0\right)\edge \left(1\ 0\ 1\right)^T$

$\left(1\ 0\ 1\right)\edge \left(1\ 0\ 0\right)^T\edge \left(1\ 0\
0\right)\edge \left(1\ 1\ 0\right)^T$ \hsk yields \hsk $\left(1\
0\ 1\right)\edge \left(1\ 1\ 0\right)^T$

$\left(1\ 0\ 1\right)\edge \left(1\ 0\ 0\right)^T\edge \left(1\ 1\
0\right)\edge \left(0\ 1\ 1\right)^T$ \hsk yields \hsk $\left(1\
0\ 1\right)\edge \left(0\ 1\ 1\right)^T$

$\left(0\ 1\ 1\right)\edge \left(1\ 1\ 0\right)^T\edge \left(1\ 0\
0\right)\edge \left(1\ 0\ 1\right)^T$ \hsk yields \hsk $\left(0\
1\ 1\right)\edge \left(1\ 0\ 1\right)^T$

$\left(0\ 1\ 0\right)\edge \left(1\ 1\ 0\right)^T\edge \left(1\ 0\
1\right)\edge \left(0\ 1\ 1\right)^T$ \hsk yields \hsk $\left(0\
1\ 0\right)\edge \left(0\ 1\ 1\right)^T$

$\left(0\ 1\ 1\right)\edge \left(1\ 0\ 1\right)^T\edge \left(1\ 1\
0\right)\edge \left(0\ 1\ 0\right)^T$ \hsk yields \hsk $\left(0\
1\ 1\right)\edge \left(0\ 1\ 0\right)^T$

$\left(0\ 1\ 0\right)\edge \left(0\ 1\ 1\right)^T\edge \left(1\ 1\
0\right)\edge \left(0\ 1\ 0\right)^T$ \hsk yields \hsk $\left(0\
1\ 0\right)\edge \left(0\ 1\ 0\right)^T$

$\left(0\ 0\ 1\right)\edge \left(1\ 0\ 1\right)^T\edge \left(1\ 1\
0\right)\edge \left(0\ 1\ 1\right)^T$ \hsk yields \hsk $\left(0\
0\ 1\right)\edge \left(0\ 1\ 1\right)^T$

$\left(0\ 0\ 1\right)\edge \left(0\ 1\ 1\right)^T\edge \left(1\ 0\
1\right)\edge \left(0\ 0\ 1\right)^T$ \hsk yields \hsk $\left(0\
0\ 1\right)\edge \left(0\ 0\ 1\right)^T$

$\left(0\ 1\ 1\right)\edge \left(1\ 1\ 0\right)^T\edge \left(1\ 0\
1\right)\edge \left(0\ 0\ 1\right)^T$ \hsk yields \hsk $\left(0\
1\ 1\right)\edge \left(0\ 0\ 1\right)^T$

\msk

For our next collection ${\mathcal V}_2$ of vertices we add
$v,w^T$ all with one entry 0, another entry 1, and the remaining
entry $a\neq 0,1$. We extend $T_1$ to a tree $T_2$ with vertex set
${\mathcal V}_2$ by adding the edge from each new vertex to the
vertex $\vei$ or $\vei^T$ (as appropriate), where $i$ is the
coordinate with entry equal to 1. In addition to the edges joining
vertices in ${\mathcal V}_1$ and those in $T_2$, the only edges
joining vertices in ${\mathcal V}_2$ are of one of the forms

\ssk

$\left(0\ 1\ a\right)\edge\left(1\ 0\ a^{-1}\right)^T$

$\left(0\ 1\ a\right)\edge\left(b\ 1\ 0\right)^T$ (with $b\neq 0$)

$\left(0\ 1\ a\right)\edge\left(0\ 1-a\ 1\right)^T$

\ssk

\noindent (together with pairs resulting from simultaneous
permutation of the coordinates of each side), since $\left(0\ 1\
a\right)\edge\left(x\ y\ z\right)^T$ requires $y+za=1$ with at
least one of $y,z$ equal to $0$ or $1$. $y=0$ implies the first
case, $y=1$ implies $z=0$ (and vice versa) and implies the second,
and $z=1$ implies the third. For these edges the 4-cycles

\ssk

\noindent $\left(0\ 1\ a\right)\edge \left(0\ 1\
0\right)^T\edge\left(0\ 1\ 0\right)\edge \left(b\ 1\ 0\right)^T$
\hsk yields \hsk $\left(0\ 1\ a\right)\edge \left(b\ 1\
0\right)^T$

\noindent together with simultaneous permutations, by permuting
throughout the 4-cycle. The 4-cycle

\ssk

\noindent $\left(0\ 1\ a\right)\edge \left(0\ 1\ 0\right)^T\edge
\left(1\ 1\ 0\right)\edge \left(1\ 0\ a^{-1}\right)^T$ \hsk yields
\hsk $\left(0\ 1\ a\right)\edge \left(1\ 0\ a^{-1}\right)^T$,

\noindent where the edge $\left(1\ 1\ 0\right)\edge \left(1\ 0\
a^{-1}\right)^T$ is a permutation of $\left(0\ 1\ 1\right)\edge
\left(a^{-1}\ 1\ 0\right)^T$ (cycling to the left), which turned
green in the previous step. We again have all simultaneous
permutations. Finally, the 4-cycle

\ssk

\noindent $\left(0\ 1\ a\right)\edge \left(1\ 1\ 0\right)^T\edge
\left(1\ 0\ 1\right)\edge \left(0\ 1-a\ 1\right)^T$ \hsk yields
\hsk $\left(0\ 1\ a\right)\edge \left(0\ 1-a\ 1\right)^T$

\noindent together with simultaneous permutations.

\msk

This deals with the full subcomplex on the vertices with at least
one 0-entry and at least one 1-entry. The next collection
${\mathcal V}_3$ of vertices adds $a\vei,a\vei^T$ with $a\neq 1$.
We extend our tree $T_2$ to a tree $T_3$ by adding the edges

$a\vei\edge\vej^T+a^{-1}\vei^T$ and $a\vei^T\edge\vej+a^{-1}\vei$

\noindent where $j$ is the smallest index $\neq i$. The edges
still unaccounted for, running between vertices of ${\mathcal
V}_3$, are

$a\vei\edge a^{-1}\vei^T$, $a\vei\edge \vek+a^{-1}\vei^T$, and
$a\vei^T\edge \vek+a^{-1}\vei$ (where $k\neq i,j$)

\noindent since $a\vei\edge x\vei^T+y\vej^T+z\vek^T$ and
$x\vei^T+y\vej^T+z\vek^T\in{\mathcal V}_3$ requires $x=a^{-1}$ and
$\{y,z\}=\{0,1\}$. The 4-cycles \ssk

$a\vei\edge\vej^T+a^{-1}\vei^T\edge\vej+\vek\edge\vek^T+a^{-1}\vei^T$,

$a\vei^T\edge\vej+a^{-1}\vei\edge\vej^T+\vek^T\edge\vek+a^{-1}\vei$,
and

$a\vei\edge\vej^T+a^{-1}\vei^T\edge\vek+a\vei\edge a^{-1}\vei^T$

\ssk

\noindent demonstrate that these edges can be turned green.

\msk

For ${\mathcal V}_4$ we add the vertices
$a\vei+b\vej,a\vei^T+b\vej^T$ with $a,b\neq 0,1$ and $i<j$. The
vertex set ${\mathcal V}_4$ thus consists of all of the vertices
with at least one 0-entry. We extend the tree $T_3$ to a tree
$T_4$ by adding the edges $a\vei+b\vej\edge a^{-1}\vei^T$ and
$a\vei^T+b\vej^T\edge a^{-1}\vei$. The edges we need to turn green
are of the form $a\vei+b\vej\edge x\vei^T+y\vej^T+z\vek^T$ and
$a\vei^T+b\vej^T\edge x\vei+y\vej+z\vek$ with $za+yb+z0=xa+yb=1$
(resp. $az+by+0z=ax+by=1$ ; we are working over a division ring!)
and at least one of $x,y,z$ equal to 0. This yields the four cases

(two coefficients equal 0): $a\vei+b\vej\edge b^{-1}\vej^T$ and
$a\vei^T+b\vej^T\edge b^{-1}\vej$,

($x=0$): $a\vei+b\vej\edge b^{-1}\vej^T+z\vek^T$ and
$a\vei^T+b\vej^T\edge b^{-1}\vej+z\vek$,

($y=0$): $a\vei+b\vej\edge a^{-1}\vei^T+z\vek^T$ and
$a\vei^T+b\vej^T\edge a^{-1}\vei+z\vek$,

($z=0$): $a\vei+b\vej\edge x\vei^T+y\vej^T$ and
$a\vei^T+b\vej^T\edge x\vei+y\vej$

\noindent These can be turned green by using the 4-cycles

[$y=0$ and $i<k$]: $a\vei+b\vej\edge a^{-1}\vei^T\edge a\vei\edge
a^{-1}\vei^T+z\vek^T$ (and transposes),

[$y=0$ and $i>k$, so $k<i<j$]: $a\vei+b\vej\edge a^{-1}\vei^T\edge
a\vei+\vej\edge z\vek^T+\vej^T$,

\hfill and so $a\vei+b\vej\edge z\vek^T+\vej^T \edge z^{-1}\vek
\edge a^{-1}\vei^T+z\vek^T$ (together with transposes),

[two coefficients 0]: $a\vei+b\vej\edge a^{-1}\vei^T+\vek^T\edge
b\vej+\vek\edge b^{-1}\vej^T$ (and transposes),

[$x=0$]: $a\vei+b\vej\edge b^{-1}\vej^T\edge b\vej\edge
b^{-1}\vej^T+z\vek^T$ (and transposes), and

[$z=0$]: $a\vei+b\vej\edge a^{-1}\vei^T+\vek^T\edge
y^{-1}\vej+\vek\edge x\vei^T+y\vej^T$ (and transposes).

\msk

Finally, ${\mathcal V}_5 = \widetilde{K}_{3,1}^{(0)}$, that is, we
add the vertices $v,w^T$ with all entries non-zero. We extend
$T_4$ to a tree $T_5$ by adding the edges

$\left(a\ b\ c\right)\edge\left(a^{-1}\ 0\ 0\right)^T$ and
$\left(a\ b\ c\right)^T\edge\left(a^{-1}\ 0\ 0\right)$.

\noindent Then the 4-cycles

$\left(a\ b\ c\right)\edge\left(a^{-1}\ 0\ 0\right)^T
\edge\left(a\ b\ 0\right)\edge\left(0\ b^{-1}\ 0\right)^T$

\hfill yield $\left(a\ b\ c\right)\edge\left(0\ b^{-1}\
0\right)^T$, and

$\left(a\ b\ c\right)\edge\left(a^{-1}\ 0\ 0\right)^T
\edge\left(a\ 0\ c\right)\edge\left(0\ 0\ c^{-1}\right)^T$

\hfill yield $\left(a\ b\ c\right)\edge\left(0\ 0\
c^{-1}\right)^T$ (together with transposes).

\noindent Then if $\left(x\ y\ z\right)$ has $z=0$ and $xa+yb=1$
with $x,y\neq 0$, the 4-cycle

$\left(a\ b\ c\right)\edge\left(a^{-1}\ 0\ 0\right)^T
\edge\left(a\ b\ 0\right)\edge\left(x\ y\ 0\right)^T$ yields
$\left(a\ b\ c\right)\edge\left(x\ y\ 0\right)^T$ (the transpose
relation is similar),

\noindent and a similar argument yields the cases where the second
or first entry is $0$. Finally, if $xa+yb+zc=1$ with $x,y,z\neq
0$, then the 4-cycle

\ssk

\ctln{$\left(a\ b\ c\right)\edge\left(a^{-1}\ 0\
0\right)^T\edge\left(a\ b+y^{-1}zc\ 0\right)\edge\left(x\ y\
z\right)^T$}

\ssk

\hfill yields $\left(a\ b\ c\right)\edge\left(x\ y\ z\right)^T$
(and the transpose relation is, again, similar).

\msk

With this, we have constructed a maximal tree $T=T_5$ in
$\widetilde{K}_{3,1}^{(1)}$, and have shown how to homotope every
edge in $\widetilde{K}_{3,1}^{(1)}$, rel endpoints, to an edge
path in $T$. Consequently, $\widetilde{K}_{3,1}$ is connected and
simply connected.

\bsk

To finish our argument, we show how to extend this result to
arbitrary $n\geq 3$. We argue by induction. The base case $n=3$ is
established above. For the inductive step, we assume that we have
shown that $\widetilde{K}_{n-1,1}$ is simply-connected. In
particular, we have constructed a maximal tree $T_{n-1}$ in
$\widetilde{K}_{n-1,1}^{(1)}$ and have shown that each of the
edges of $\widetilde{K}_{n-1,1}^{(1)}$ not in $T_{n-1}$ can be
deformed, rel endpoints, in $\widetilde{K}_{n-1,1}$, into
$T_{n-1}$.

By appending $0$'s to the end of every column of the vectors $v,w$
labelling the vertices $v,w^T$ of $\widetilde{K}_{n-1,1}$
(yielding matrices $v_+,w_+$) and noting that $w_+^Tv_+=1$ iff
$w^Tv=1$, the map $v\mapsto v_+$, $w^T\mapsto w_+^T$ induces an
embedding of $\widetilde{K}_{n-1,1}$ into $\widetilde{K}_{n,1}$,
and its image is the full subcomplex of $\widetilde{K}_{n,1}$ on
the vertex set ${\mathcal V}_6$ = the image of
$\widetilde{K}_{n-1,1}^{(0)}$. The image of the tree $T$ is a tree
$T_6$ which provides the starting point for constructing the
needed tree in $\widetilde{K}_{n,1}$.

To build our tree $T_7$ we add, for the vertices
$\left(\begin{matrix} a_1&\ldots&a_n
\end{matrix}\right),\left(\begin{matrix} a_1&\ldots&a_n
\end{matrix}\right)^T$ with $a_n\neq 0$ and at least one other entry
$a_j\neq 0$ (we may assume $j$ is the smallest such index), the
edges $\left(\begin{matrix} a_1&\ldots&a_n
\end{matrix}\right)\edge a_j^{-1}\vej^T$ and $\left(\begin{matrix}
a_1&\ldots&a_n
\end{matrix}\right)^T\edge a_j^{-1}\vej$, and then add the edges
$a_n\ven\edge\veu^T+a_n^{-1}\ven^T$ and
$a_n\ven^T\edge\veu+a_n^{-1}\ven$ for each $a_n\neq 0$. This gives
a maximal tree in $\widetilde{K}_{n,1}$; by our inductive step we
know that every edge in the image of $\widetilde{K}_{n-1,1}^{(1)}$
is homotopic, rel endpoints, to an edge path in $T_7$, and so can
be assumed to be green.

We now work our way through all of the remaining edges of
$\widetilde{K}_{n,1}$ in steps, to show that they can all be
turned green. If a vertex $v,w^T$ has two or more entries
$a_j,a_k\neq 0$, for $j,k<n$ (we may assume $j$ is the smallest
such index), then the 4-cycles

\ssk

\ctln{$v\edge a_j^{-1}\vej^T\edge a_j\vej+a_k\vek\edge
a_k^{-1}\vek^T$ and $w^T\edge a_j^{-1}\vej\edge
a_j\vej^T+a_k\vek^T\edge a_k^{-1}\vek$}

\ssk

\noindent enable us to make the edges $v\edge a_k^{-1}\vek^T$ and
$w^T a_k^{-1}\vek$ green, and allowing us to base our further
arguments off of any non-zero entry of $v,w^T$ (other than the
last entry). If $v=\left(\begin{matrix} a_1&\ldots&a_n
\end{matrix}\right)\edge\left(\begin{matrix} x_1&\ldots&x_{n-1}&0
\end{matrix}\right)^T=w^T$ is an edge, then $x_ia_i\neq 0$ for some
$i$, and the 4-cycle

\ssk

\ctln{$v\edge a_i^{-1}\vei^T\edge \left(\begin{matrix}
a_1&\ldots&a_{n-1}&0 \end{matrix}\right)\edge w^T$}

\ssk

\noindent enables us to turn the edge $v\edge w^T$ green, and a
similar argument will allow us to turn the edges
$v=\left(\begin{matrix} x_1&\ldots&x_{n-1}&0
\end{matrix}\right)\edge\left(\begin{matrix} a_1&\ldots&a_n
\end{matrix}\right)^T=w^T$ green. With this we can make every edge
of $\widetilde{K}_{n,1}$ joining a vertex of
$\widetilde{K}_{n,1}$, \underbar{other} \underbar{than} the
vertices $a_n\ven,a_n\ven^T$, to a vertex in the image of
$\widetilde{K}_{n-1,1}^{(0)}$, green. Note, however, that there
are no edges between the vertices $a_n\ven,a_n\ven^T$ and the
vertices in the image of $\widetilde{K}_{n-1,1}^{(0)}$.

\msk

For every $i>1$ the 4-cycle $a_n\ven\edge
\veu^T+a_n^{-1}\ven^T\edge\veu+\vei\edge\vei^T+a_n^{-1}\ven^T$
makes $a_n\ven\edge\vei^T+a_n^{-1}\ven^T$ green, and a similar
argument turns $a_n\ven^T\edge\vei+a_n^{-1}\ven$ green. Then for
$i\neq j$, $i,j<n$,

\ssk

\ctln{$\vei+a_n\ven\edge
\vei^T+\vej^T\edge\vej\edge\vej^T+a_n^{-1}\ven^T$}

\ssk

\noindent turns $\vei+a_n\ven\edge\vej^T+a_n^{-1}\ven^T$ green.
Then (using the fact that $n\geq 3$)

\ssk

\ctln{$a_n\ven\edge \veu^T+a_n^{-1}\ven^T\edge\vet+a_n\ven\edge
a_n^{-1}\ven^T$}

\ssk

\noindent turns $a_n\ven\edge a_n^{-1}\ven^T$ green. If one of
$x_1,\ldots x_{n-1}$ is non-zero (say $x_i$, and then choose
$j\neq i$, $j<n$), then the 4-cycles

\ssk

\ctln{$a_n\ven\edge \vej^T+a_n^{-1}\ven^T\edge
x_i^{-1}(1-x_j)\vei+\vej\edge\left(\begin{matrix}
x_1&\ldots&x_{n-1}&a_n^{-1} \end{matrix}\right)^T$ and}

\ssk

\ctln{$a_n\ven^T\edge \vej+a_n^{-1}\ven\edge
x_i^{-1}(1-x_j)\vei^T+\vej^T\edge\left(\begin{matrix}
x_1&\ldots&x_{n-1}&a_n^{-1} \end{matrix}\right)$}

\ssk

\noindent together with the previous step enable us to make every
edge with vertex either $a_n\ven$ or $a_n\ven^T$ green.

\msk

All that is left now is to turn the edges $v=\left(\begin{matrix}
a_1&\ldots&a_n \end{matrix}\right)\edge\left(\begin{matrix}
x_1&\ldots&x_n
\end{matrix}\right)^T=w^T$ with $a_n,x_n\neq 0$ green. By the
immediately preceding step we may also assume that $x_i,a_j\neq 0$
for some $i,j<n$. If $i\neq j$ then the 4-cycle

\ssk

\ctln{$v\edge a_j^{-1}\vej^T\edge a_j\vej+x_i^{-1}(1-a_jx_j)\vei\edge w^T$}

\ssk

\noindent turns the edge $v\edge w^T$ green. If we cannot find
such a pair of distinct indices $i,j$ then our edge must be of the
form $a_i\vei+a_n\ven\edge x_i\vei^T+x_n\ven^T$ with
$x_ia_i+x_na_n=1$ and $a_i,x_i,a_n,x_n\neq 0$. Then (choosing a
$j\neq i,n$) the 4-cycle

\ssk

\ctln{$a_i\vei+a_n\ven\edge a_i^{-1}\vei^T+\vej^T\edge
x_i^{-1}\vei+(1-a_i^{-1}x_i^{-1})\vej\edge x_i\vei^T+x_n\ven^T$}

\ssk

\noindent demonstrates that this final collection of edges can be
made green.

\bsk

Therefore, every edge of $\widetilde{K}_{n,1}^{(1)}$ is homotopic,
rel endpoints, to an edge path in $T_7$, and $\widetilde{K}_{n,1}$
is simply connected. This finishes the inductive step; so for
every $n\geq 3$, $\widetilde{K}_{n,1}$ is connected and simply
connected, and so is the universal covering space of $K_{n,1}$.
This completes the proof of Theorem \ref{rankone}.

\bsk

\section{Closing Remarks}

Theorem \ref{rankone} provides a natural example of a torsion
group that arises as a maximal subgroup of the free idempotent
generated semigroup on some (finite) biordered set, answering a
question raised in \cite{ESV}. After the results of this paper were
announced, Gray and Ruskuc \cite{GR} proved that every group arises 
as the maximal subgroup of some biordered set overriding this particular 
example. We remark that if $e$ is an
idempotent matrix of rank $n-1$ in $E = E(M_{n}(Q))$, then the
maximal subgroup of $IG(E)$ with identity $e$ must be a free group
by Theorem \ref{BMM}, since there are no idempotents available to
singularize a square consisting of rank $n-1$ idempotent matrices.
Thus the maximal subgroup of $IG(E)$ corresponding to an
idempotent of rank $n-1$ is not isomorphic to $GL_{n-1}(Q)$.

Based on experimental evidence, we conjecture that the maximal
subgroup of $IG(E)$ with identity an idempotent matrix of rank $k
< n-1$ is $GL_{k}(Q)$, at least if $k < n/2$ and $n \geq 3$, but
this problem remains open. It is plausible that when $k < n-1$ the
subcomplex $\widetilde{K}_{n,k}$ of $K$ spanned by the vertices of
rank $k$ is simply connected. The methods in the proof of Theorem
\ref{rankone} in the present paper seem difficult to extend.
However, there is a lot of structure to the complexes that we have
not exploited. The connections to  Grassmanians - the vertices of
${K}_{n,k}$ are (two sets of) the points of the Grassmanian
$G_{n,k}$ of $k$-planes in $Q^n$, and the vertices of
$\widetilde{K}_{n,k}$ are (two sets of) the points of the
universal bundle over $G_{n,k}$ - seem worth exploring further,
and we expect to consider these ideas in a subsequent paper.

\end{document}